\documentclass{article} \UseRawInputEncoding
\usepackage{amsmath,amsthm,amssymb,amsfonts}
\usepackage{verbatim}
\usepackage{graphicx}
\usepackage{stmaryrd} 
\usepackage{hyperref}

\usepackage[colorinlistoftodos]{todonotes}

\theoremstyle{plain}
\newtheorem{thm}{Theorem}
\newtheorem{lem}[thm]{Lemma}
\newtheorem{prop}[thm]{Proposition}
\newtheorem{cor}[thm]{Corollary}

\theoremstyle{definition}
\newtheorem{definition}[thm]{Definition}
\newtheorem{exl}[thm]{Example}

\numberwithin{thm}{section}

\newcommand{\adj}{\leftrightarrow}
\newcommand{\adjeq}{\leftrightarroweq}

\DeclareMathOperator{\id}{id}
\def\Z{{\mathbb Z}}
\def\N{{\mathbb N}}

\def \Fix {Fix}

\title{Irreducibility 
and Rigidity in 
Digital Images}

\author{Laurence Boxer
\thanks{Department of Computer and Information Sciences, Niagara University, NY 14109, USA
and  \newline
Department of Computer Science and Engineering, State University of New York at Buffalo \newline
email: boxer@niagara.edu
}
}

\date{ }
\begin{document}
\maketitle{}

\begin{abstract}
We study how
the properties of
irreducibility and rigidity
in digital images interact
with Cartesian products,
wedges, and cold and
freezing sets.
\newline

Key words and phrases: digital topology, digital image, freezing set,
reducible, rigid
\newline

MSC: 54B20, 54C35
\end{abstract}

\date
\maketitle

\section{Introduction}
The properties
of irreducibility
and rigidity in
digital images
were introduced
in~\cite{hmps} and
have been studied
in subsequent papers,
including~\cite{bs19a,BxFPSets2,BxColdFreeze,BxLtd}.
In the current work, we
study implications of
these properties for
Cartesian products,
wedges, and cold and
freezing sets.

\section{Preliminaries}
We use $\N$ for the set of natural numbers,
$\N^* = \N \cup \{0\}$,
$\Z$ for the set of integers, and
$\#X$ for the number of distinct members of $X$.

We typically denote a (binary) digital image
as $(X,\kappa)$, where $X \subset \Z^n$ for some
$n \in \N$ and $\kappa$ represents an adjacency
relation of pairs of points in $X$. Thus,
$(X,\kappa)$ is a graph, in which members of $X$ may be
thought of as black points, and members of $\Z^n \setminus X$
as white points, of a picture of some ``real world" 
object or scene.

\subsection{Adjacencies}
This section is largely
quoted or paraphrased
from~\cite{BxArbDim}.

Let $u,n \in \N$, $1 \le u \le n$. 
For $X \subset \Z^n$,
$x = (x_1, \ldots, x_n),~y=(y_1, \ldots, y_n) \in X$ 
are $c_u$-adjacent if and only if
\begin{itemize}
    \item $x \neq y$, and
    \item for at most $u$ indices~$i$, 
          $\mid x_i - y_i \mid = 1$, and
    \item for all indices $j$ such that 
          $\mid x_j - y_j \mid \neq 1$, we have
          $x_j = y_j$.
\end{itemize}
The $c_u$ adjacencies are the adjacencies most used
in digital topology, especially $c_1$ and $c_n$.

In low dimensions, it is also common to denote a
$c_u$ adjacency by the number of points that can
have this adjacency with a given point in $\Z^n$. E.g.,
\begin{itemize}
    \item 
    in $\Z$, $c_1$-adjacency is 2-adjacency;
    \item in 
    $\Z^2$, $c_1$-adjacency is 4-adjacency and
$c_2$-adjacency is 8-adjacency;
\item  in 
$\Z^3$, $c_1$-adjacency is 6-adjacency,
    $c_2$-adjacency is 18-adjacency, and
    $c_3$-adjacency is 26-adjacency.
\end{itemize}

We use the notations $y \adj_{\kappa} x$, or, when
the adjacency $\kappa$ can be assumed, $y \adj x$, to mean
$x$ and $y$ are $\kappa$-adjacent.
The notations $y \adjeq_{\kappa} x$, or, when
$\kappa$ can be assumed, $y \adjeq x$, mean either
$y=x$ or $y \adj_{\kappa} x$.

A sequence $P=\{y_i\}_{i=0}^m$ in a digital image $(X,\kappa)$ is
a {\em $\kappa$-path from $a \in X$ to $b \in X$} if
$a=y_0$, $b=y_m$, and $y_i \adjeq_{\kappa} y_{i+1}$ 
for $0 \leq i < m$.

$X$ is {\em $\kappa$-connected}~\cite{Rosenfeld},
or {\em connected} when $\kappa$
is understood, if for every pair of points $a,b \in X$ there
exists a $\kappa$-path in $X$ from $a$ to $b$.

A {\em (digital) $\kappa$-closed curve} is a
path $S=\{s_i\}_{i=0}^{m-1}$ such that $s_0 \adj_{\kappa} s_{m-1}$,
and $i \neq j$ 
implies $s_i \neq s_j$. If also $0 \le i < m$ implies
the only $\kappa$-adjacent members of $S$ to $x_i$ are 
$x_{(i-1) \mod m}$ and $x_{(i+1) \mod m}$,
then $S$ is a {\em (digital) 
$\kappa$-simple closed curve}. 

\subsection{Digitally continuous functions}
This section is largely
quoted or paraphrased
from~\cite{BxArbDim}.

Digital continuity is defined
to preserve connectedness, as at
Definition~\ref{continuous} below. By
using adjacency as our standard of ``closeness," we
get Theorem~\ref{continuityPreserveAdj} below.

\begin{definition}
\label{continuous}
{\rm ~\cite{Bx99} (generalizing a definition of~\cite{Rosenfeld})}
Let $(X,\kappa)$ and $(Y,\lambda)$ be digital images.
A function $f: X \rightarrow Y$ is 
{\em $(\kappa,\lambda)$-continuous} if for
every $\kappa$-connected $A \subset X$ we have that
$f(A)$ is a $\lambda$-connected subset of $Y$.
\end{definition}

When $X \cup Y \subset
(Z^n, \kappa)$, we use the abbreviation
{\em $\kappa$-continuous} for {\em 
$(\kappa,\kappa)$-continuous}.

When the adjacency relations are understood, we will simply say that $f$ is \emph{continuous}. Continuity can be expressed in terms of adjacency of points:

\begin{thm}
{\rm ~\cite{Rosenfeld,Bx99}}
\label{continuityPreserveAdj}
A function $f:X\to Y$ is continuous if and only if $x \adj x'$ in $X$ 
implies $f(x) \adjeq f(x')$.
\end{thm}

See also~\cite{Chen94,Chen04}, where similar notions are referred to as {\em immersions}, {\em gradually varied operators},
and {\em gradually varied mappings}.

A digital {\em isomorphism} (called {\em homeomorphism}
in~\cite{Bx94}) is a $(\kappa,\lambda)$-continuous
surjection $f: X \to Y$ such that $f^{-1}: Y \to X$ is
$(\lambda,\kappa)$-continuous.

A {\em digital interval} is
a set denoted $[a,b]_{\Z}$
where $a,b \in \Z$, 
$a \le b$, and
\[ [a,b]_{\Z} = \{z \in \Z
 \, | \, a \le z \le b \}
\]
with the $c_1$ adjacency
in~$\Z$.

Let $X \subset \Z^n$.
The
{\em boundary of} $X$
{\rm \cite{Rosenf79}} is
\[Bd(X) = \{x \in X \, | \mbox{ there exists } y \in \Z^n \setminus X \mbox{ such that } y \adj_{c_1} x\}.
\]

A homotopy between continuous functions may be thought of as
a continuous deformation of one of the functions into the 
other over a finite time period.

\begin{definition}{\rm (\cite{Bx99}; see also \cite{Khalimsky})}
\label{htpy-2nd-def}
Let $X$ and $Y$ be digital images.
Let $f,g: X \rightarrow Y$ be $(\kappa,\kappa')$-continuous functions.
Suppose there is a positive integer $m$ and a function
$F: X \times [0,m]_{{\Z}} \rightarrow Y$
such that

\begin{itemize}
\item for all $x \in X$, $F(x,0) = f(x)$ and $F(x,m) = g(x)$;
\item for all $x \in X$, the induced function
      $F_x: [0,m]_{{\Z}} \rightarrow Y$ defined by
          \[ F_x(t) ~=~ F(x,t) \mbox{ for all } t \in [0,m]_{{\Z}} \]
          is $(2,\kappa')-$continuous. That is, $F_x(t)$ is a path in $Y$.
\item for all $t \in [0,m]_{{\Z}}$, the induced function
         $F_t: X \rightarrow Y$ defined by
          \[ F_t(x) ~=~ F(x,t) \mbox{ for all } x \in  X \]
          is $(\kappa,\kappa')-$continuous.
\end{itemize}
Then $F$ is a {\rm digital $(\kappa,\kappa')-$homotopy between} $f$ and
$g$, and $f$ and $g$ are {\rm digitally $(\kappa,\kappa')-$homotopic in} $Y$.
$\Box$
\end{definition}

\begin{thm}
    \label{SCChtpyPreservesOnto}
    {\rm \cite{Bx10}}
     Let $S$ be a simple closed 
    $\kappa$-curve
and let $H: S \times [0,m]_{\Z}
\to S$ be a 
$(\kappa,\kappa)$-homotopy 
between an isomorphism $H_0$ and
$H_m = f$, where
$f(S) \neq S$. Then $\#S = 4$.
\end{thm}

The literature uses {\em path} polymorphically: a
$(c_1,\kappa)$-continuous 
function $f: [0,m]_{\Z} \to X$
is a $\kappa$-path if 
$f([0,m]_{\Z})$ is a 
$\kappa$-path from $f(0)$ 
to $f(m)$
as described above.

We use $\id_X$ to denote
the {\em identity function}, $\id_X(x) = x$
for all $x \in X$.

Given a digital image
$(X,\kappa)$, we denote
by $C(X,\kappa)$ the set
of $\kappa$-continuous
functions $f: X \to X$.

Given $f \in C(X,\kappa)$,
a {\em fixed point} of $f$
is a point $x \in X$ such
that $f(x)=x$.
$\Fix(f)$ will denote
the set of fixed points
of~$f$. We say
$f$ is a {\em retraction},
and the set $Y=f(X)$ is a
{\em retract of $X$}, if
$f|_Y = \id_Y$; thus, 
$Y = \Fix(f)$.

\begin{definition}
\label{freezeDef}
{\rm \cite{BxFPSets2}}
Let $(X,\kappa)$ be a
digital image. We say
$A \subset X$ is a 
{\em freezing set for $X$}
if given $g \in C(X,\kappa)$,
$A \subset \Fix(g)$ implies
$g=\id_X$. A freezing set
$A$ is {\em minimal} if
no proper subset of $A$
is a freezing set for
$(X,\kappa)$.
\end{definition}

\begin{exl} We have
the following examples
from~{\rm \cite{BxFPSets2}}.
\begin{itemize}
    \item  $\{a,b\}$ is a minimal
    freezing set for
    $[a,b]_{\Z}$.
    \item
    Given $X \subset \Z^n$
    such that $X$ is 
    finite and 
    $1 \le u \le n$,
    $Bd(X)$ is a freezing
    set for $(X,c_u)$
    (not necessarily minimal).
    \item $\Pi_{i=1}^n \{a_i,b_i \}$ is a freezing set for
    $(X,c_1)$, where
    $X = \Pi_{i=1}^n [a_i,b_i]_{\Z}$ 
    (minimal for $n \in \{1,2\}$; not 
    necessarily minimal for
    $n>2$).
 \end{itemize}   
\end{exl}

The following elementary
assertion was noted
in~\cite{BxFPSets2}.

\begin{lem}
\label{swelling}
Let $(X,\kappa)$ be
    a connected digital 
    image for which $A$ is
    a freezing set. If
    $A \subset A' \subset X$, then $A'$ is a
    freezing set for 
    $(X,\kappa)$.
\end{lem}

\begin{definition}
{\rm \cite{BxFPSets2}}
\label{s-cold-def}
Given $s \in \N^*$, we say $A \subset X$ is an
{\em $s$-cold set} for the connected digital image $(X,\kappa)$
if given $g \in C(X,\kappa)$ such that
$g|_A = \id_A$, then for all $x \in X$, there is a
$\kappa$-path in~$X$ of
length at most~$s$ from~$x$
to~$g(x)$.
A {\em cold set} is a 1-cold set.
\end{definition}

\begin{exl}
    {\rm \cite{BxFPSets2}}
    $\{0\}$ is a cold set,
    but not a freezing set,
    for $[0,1]_{\Z}$.
\end{exl}

Note a 0-cold set is a freezing set~{\rm \cite{BxFPSets2}}.

Let $X \subset \Z^n$,
$x = (x_1, \ldots, x_n) \in Z^n$, where
each $x_i \in \Z$. For
each index~$i$,
the {\em projection map} (onto the $i^{th}$
coordinate) $p_i: X \to \Z$ is given by
$p_i(x) = x_i$.

\subsection{Tools for
determining fixed point sets}
\begin{thm}
\label{freezeInvariant}
{\rm \cite{BxFPSets2}}
Let $A$ be a freezing set for the digital image $(X,\kappa)$ and let
$F: (X,\kappa) \to (Y,\lambda)$ be an isomorphism. Then $F(A)$ is
a freezing set for $(Y,\lambda)$.
\end{thm}

\begin{prop}
{\rm \cite{bs19a}}
\label{uniqueShortest}
Let $(X,\kappa)$ be a digital
image and $f \in C(X,\kappa)$.
Suppose $x,x' \in \Fix(f)$ are
such that there is a unique
shortest $\kappa$-path $P$ in $X$ 
from $x$ to $x'$. Then
$P \subset \Fix(f)$.
\end{prop}

The following lemma may be
understood as saying that
if $q$ and $q'$ are
adjacent with $q$ in a
given direction from $q'$,
and if $f$ pulls $q$ 
further in that direction,
then $f$ also pulls $q'$
in that direction.

\begin{lem}
\label{pullingLem}
{\rm \cite{BxFPSets2}}
Let $(X,c_u)\subset \Z^n$ be a digital image, 
$1 \le u \le n$. Let $q, q' \in X$ be such that
$q \adj_{c_u} q'$.
Let $f \in C(X,c_u)$.
\begin{enumerate}
    \item If $p_i(f(q)) < p_i(q) < p_i(q')$
          then $p_i(f(q')) < p_i(q')$.
    \item If $p_i(f(q)) > p_i(q) > p_i(q')$
          then $p_i(f(q')) > p_i(q')$.
\end{enumerate}
\end{lem}

\subsection{Irreducible and Rigid Images}
\begin{definition}
{\rm \cite{hmps}}
A finite image $X$ is {\em reducible} when it is
homotopy equivalent to an image of fewer points. Otherwise,
we say $X$ is 
{\em irreducible}.
\end{definition}

\begin{lem}
    {\rm \cite{hmps}}
    \label{reducibleCharacterize}
    A finite image $X$ is
    reducible if and only 
    if $\id_X$ is 
    homotopic to a 
    nonsurjective map.
\end{lem}

\begin{figure}[htb]
\includegraphics[width=99mm]{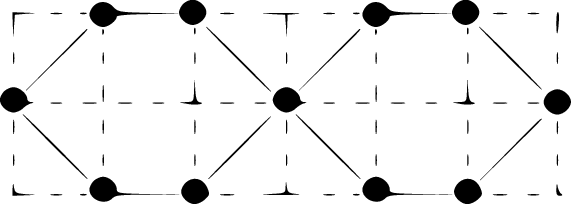}
\caption{\cite{BxSt19}
    Example of a rigid
    digital image - a wedge of
    digital simple closed curves}
\label{fig:rigid}
\end{figure}

\begin{lem}
{\rm \cite{hmps}}
\label{1step}
    A finite image $X$ is 
    reducible if and only if
    $\id_X$ is homotopic in 
    one step to a 
    nonsurjective map.
\end{lem}

\begin{definition}
\label{rigidDef}
{\rm \cite{hmps}}
    We say an image $X$
    is {\em rigid} if the 
    only map
homotopic to $id_X$ 
is $id_X$.
\end{definition}

Figure~\ref{fig:rigid} shows an
example of a rigid digital image.

\begin{prop}
\label{rigidIsIrreduc}
    {\rm \cite{hmps}} A finite 
    rigid digital image is 
irreducible.
\end{prop}

That the converse of 
Proposition~\ref{rigidIsIrreduc}
is not generally valid, is shown
by the following example.

\begin{exl}
    {\rm \cite{hmps}} A
    digital simple closed curve is
    irreducible but not rigid.
\end{exl}

\section{Products}
For Cartesian products of digital images 
$X=\Pi_{i=1}^v (X_i,\kappa_i)$
and $1 \le u \le v$, we often 
use the generalized normal 
product adjacency~\cite{BxGen} 
$NP_u(\kappa_1, \ldots, \kappa_v)$: given distinct $x,x' \in X$, 
$x = (x_1, \ldots, x_v)$,
$x'=(x_1', \ldots, x_v')$, 
where $x_i, x_i' \in X_i$, we have
$x \adj_{NP_u(\kappa_1, \ldots, \kappa_v)} x'$
if and only if
\begin{itemize}
\item for at least 1 and at most $u$ indices $i$,
$x_i \adj_{\kappa_i} x_i'$, and
\item for all other indices $j$, $x_j = x_j'$.
\end{itemize}

\begin{thm}
{\rm \cite{bs19a}}
\label{reducibleProd}
    Let $(X_i,\kappa_i)$
be a digital image,
$1 \le i \le v$.
Let $X = \Pi_{i=1}^v X_i$.
If $(X,NP_v(\kappa_1,\ldots,\kappa_v))$ is rigid,
then each $(X_i,\kappa_i)$
is rigid.
\end{thm}

At Corollary~\ref{prodIrred}
below, we
obtain an analogous result
for irreducible digital images.

\begin{thm}
\label{reduProdImpliesReduFactor}
Let $(X_i,\kappa_i)$
be a finite digital image,
$1 \le i \le v$.
Let $X = \Pi_{i=1}^v X_i$.
If for some~$j$,
$(X_j,\kappa_j)$ is
reducible, then
$(X,NP_v(\kappa_1,\ldots,\kappa_v))$ is reducible.
\end{thm}

\begin{proof}
    By Lemma~\ref{reducibleCharacterize}, there is
    a $\kappa_j$-homotopy
    $H_j: X_j \times [0,m]_{\Z} \to X_j$
    from $\id_{X_j}$ to
    a nonsurjective map
    $f_j: X_j \to X_j$.
    For $i \neq j$, let
    $H_i: X_i \times [0,m]_{\Z} \to X_i$
    be the trivial
    homotopy $H_i(x_i,t) = x_i$. Then
    $H: X \times [0,m]_{\Z} \to X$, 
    given by
    \[ H(x_1, \ldots, x_v, t) = 
    (H_1(x_1,t), \ldots, H_v(x_v,t))\]
    is an $NP_v(\kappa_1,\ldots,\kappa_v)$-homotopy 
    from $\id_X$ to a
    nonsurjective map.
    The assertion follows
    from Lemma~\ref{reducibleCharacterize}.
\end{proof}

\begin{exl}
    Let $(X_1,c_2)$ be the
    rigid digital image of
    Figure~\ref{fig:rigid}.
    By Proposition~\ref{rigidIsIrreduc}, $(X_1,c_2)$
    is irreducible. Let
    $(X_2,c_1)=[0,1]_{\Z}$.
    Clearly, $(X_2,c_1)$ is
    reducible. By
    Theorem~\ref{reduProdImpliesReduFactor},
    $(X_1 \times X_2,
    NP_2(c_2,c_1))$ is
    reducible.
\end{exl}

As an immediate consequence of
Theorem~\ref{reduProdImpliesReduFactor}, we have the following.

\begin{cor}
\label{prodIrred}
    Let $(X_i,\kappa_i)$
be a finite digital image,
$1 \le i \le v$.
Let $X = \Pi_{i=1}^v X_i$. If
$(X,NP_v(\kappa_1,\ldots,\kappa_v))$ is irreducible, then each
$(X_i,\kappa_i)$ is irreducible.
\end{cor}

\section{Wedges}
Let $X \cup Y \subset (\Z^n,\kappa)$
such that there
is a point
$x_0 \in \Z^n$
with $X \cap Y = \{x_0\}$. Suppose
$x \in X$, $y \in Y$, and
$x \adjeq_{\kappa} y$ imply
$x_0 \in \{x,y\}$. Then
$X' = X \cup Y$
is the $(\kappa-)wedge$ of $X$ and $Y$,
denoted $X' = X \vee Y$. We call
$x_0$ the
{\em wedge point}
of $X'$.

In this section,
we explore the
preservation of
irreducibility and
of rigidity by
the wedge construction.

\begin{lem}
\label{wedgeRetraction}    
 Let 
    $(X,\kappa) =
    (X_0,\kappa) \vee (X_1,\kappa)$
    where $x_0$ is
    the wedge point.
    The function
    $r: X \to X_0$
    given by
    \[ r(x) = \left \{
    \begin{array}{ll}
       x  & \mbox{if } x \in X_0; \\
       x_0  & \mbox{if } x \not \in X_0,
    \end{array}
    \right .
    \]
    is $\kappa$-continuous and
    is a $\kappa$-retraction.
\end{lem}

\begin{proof}
    Elementary and left to
    the reader.
\end{proof}

We have the following.

\begin{thm}
\label{rigidWedge}
{\rm \cite{bs19a}}
        Let 
    $(X,\kappa) =
    (X_0,\kappa) \vee (X_1,\kappa)$
    where $x_0$ is
    the wedge point. 
    Suppose
    $\#X_0 > 1$ and
    $\#X_1 > 1$.
    Suppose $(X_0,\kappa)$
    and $(X_1,\kappa)$ are
    both connected.
   If $X_0$ and $X_1$ are
    both rigid, then 
    $X$ is rigid.
\end{thm}

We obtain a similar result
for the property of
irreducibility in the 
following.
    
\begin{thm}
\label{irreducibleWedge}
    Let 
    $(X,\kappa) =
    (X_0,\kappa) \vee (X_1,\kappa)$
    where $x_0$ is
    the wedge point,
    i.e., 
    $\{x_0\} = X_0 \cap X_1$. Suppose
    $\#X_0 > 1$ and
    $\#X_1 > 1$.
    If $X_0$ and $X_1$ are
    both irreducible, then 
    $X$ is irreducible.
\end{thm}

\begin{proof}
Suppose otherwise. Then there
is a digital homotopy
\[ H: X \times [0,m]_{\Z} 
\to X
\]
between $\id_X$ 
and a continuous function
$f: X \to X$ such that
$f$ is not a surjection.
Without loss of generality,
there exists $y \in X_0$
such that $y \not \in f(X)$.

Let $R$ be the retraction
of Lemma~\ref{wedgeRetraction}. Then
$R \circ H: X_0 \times [0,m]_{\Z} 
\to X_0$ is a $\kappa$-homotopy from $\id_{X_0}$
to $R \circ f|_{X_0}$, and
$y \not \in R \circ f(X_0)$.
By Lemma~\ref{reducibleCharacterize}, this is
contrary to the 
assumption that $X_0$ is 
irreducible. The assertion
follows.
\end{proof}

The converse of 
Theorem~\ref{rigidWedge}
is not generally valid, as
shown by Example~3.11
of~\cite{bs19a}.

\begin{prop}
    (Corollary~3.13 
    of~{\rm \cite{hmps}})
    \label{irredButNotRigid}
    A digital simple closed curve
    of at least 5 points is
    irreducible but not rigid.
\end{prop}

For the following
Example~\ref{irrWedgeSCC}
and Theorem~\ref{rigidWedgeSCC},
we have
\begin{itemize}
    \item $(X,\kappa)=(Y,\kappa) \vee (S,\kappa)$, where
    $\#Y > 1$,
    $(Y,\kappa)$ is 
    irreducible or rigid, and
    $(S,\kappa)$ is a digital
    simple closed curve 
    of at least 5 points.
    \item $S = \{s_i\}_{i=0}^n$
    is a circular
    listing of the 
    members of $S$,
    where $s_0=x_0$.
    \item Functions $R,~R_1: X \to X$ are given by
    \[ R(x) = \left \{ \begin{array}{ll}
        x_0 & \mbox{if } x \in Y; \\
        x & \mbox{if } x \in S,
    \end{array} \right \},
    ~~~~~~~~~~
     R_1(x) = \left \{ \begin{array}{ll}
        x_0 & \mbox{if } x \in S; \\
        x & \mbox{if } x \in Y
    \end{array} \right \}.
    \]
    \item Given a homotopy
    $H: X \times [0,m]_{\Z} \to X$ from $\id_X$ to
    $f \in C(X,\kappa)$,
    let $G: S \times [0,1]_{\Z} \to S$
be given by 
\[ G(s,t) = R(H(s,t))
\]
and let
$G_1: Y \times [0,1]_{\Z} \to Y$
be given by 
\[ G_1(x,t) = R_1(H(x,t)).
\]
\end{itemize}

\begin{exl}
\label{irrWedgeSCC}
    Let $(X,\kappa)=(Y,\kappa) \vee (S,\kappa)$, where
    $\#Y > 1$,
    $(Y,\kappa)$ is irreducible, and
    $(S,\kappa)$ is a digital
    simple closed curve 
    of at least 5 points. Then
    $(X,\kappa)$ is irreducible.
\end{exl}

\begin{proof}
The assertion follows from
Theorem~\ref{irreducibleWedge}
and Proposition~\ref{irredButNotRigid}.
\end{proof}

\begin{thm}
\label{rigidWedgeSCC}
    Let $(X,\kappa)=(Y,\kappa) \vee (S,\kappa)$, where
    $Y$ is finite and
    $\#Y > 1$,
    $(Y,\kappa)$ is rigid, and
    $(S,\kappa)$ is a digital
    simple closed curve 
    of at least 5 points. Then
    $(X,\kappa)$ is rigid.
\end{thm}

\begin{proof}
    We argue by contradiction.
    Suppose $f \in C(X,\kappa)$
    such that $f \neq \id_X$ and
    there is a homotopy
    $H: X \times [0,m]_{\Z} \to X$
    from $\id_X$ to $f$. By
    Definition~\ref{rigidDef},
    we may assume $m=1$.
    
    Let $x_0$ be the wedge point, i.e., $\{x_0\} = Y \cap S$,
    where $\{x_i\}_{i=0}^{n-1}$ is a
    circular ordering of the 
    distinct members of $S$.
    Consider the following cases.
\begin{itemize}
    \item $f(x_0)=H(x_0,1) \in Y \setminus \{x_0\}$.
    Then we must have 
    $H(x_1,1)=x_0$ and 
    $H(x_{n-1},1)=x_0$. 
    
By Lemma~\ref{wedgeRetraction},
$R$ is a retraction of~$X$ to~$S$.
We have
\begin{equation}
\label{nonOnto}
    R(f(x_0)) = x_0 = R(f(x_1))
\end{equation}
    
Then $G$ is a homotopy from $\id_S$
to a map that is non-injective, 
hence non-surjective; this is
impossible by 
Proposition~\ref{irredButNotRigid} and
Lemma~\ref{reducibleCharacterize}.

\item $f(x_0)=H(x_0,1) \in S \setminus
\{x_0\}$.

By Lemma~\ref{wedgeRetraction},
$R_1$ is a retraction. Since
$Y$ is connected and has more than
$1$ point, there exists $y \in Y$
such that $y \adj x_0$. However, $y$
is not adjacent to
any member of $S$
other than $x_0$.
Therefore,
$H(y,1) \in Y$.
Hence
\[ x_0 \adj H(x_0,1) \adj H(y,1) = x_0 \]
and 
\begin{equation}
\label{notId1}
  G_1(x_0,1) =  R_1(H(x_0,1)) = x_0 = R_1(H(y,1))
\end{equation}

Then $G_1$
is a homotopy from $\id_Y$ to a
map that, by~(\ref{notId1}), 
is not $\id_Y$. This is impossible,
since $Y$ is rigid.

\item $f(s)=H(s,1) \in Y 
\setminus \{x_0\}$ for some 
$s \in S \setminus \{x_0\}$.
This is impossible, as the only member of
$S$ that is within 1 step of
$Y \setminus \{x_0\}$ is $x_0$.
\item $f(y)=H(y,1) \in S \setminus \{x_0\}$ for some 
$y \in Y \setminus \{x_0\}$. 
This is impossible, as the only 
member of $Y$ that is within~1 
step of $S\setminus \{x_0\}$ is $x_0$.
\item $f(x_i)=H(x_i,1) =x_j$ for some 
indices satisfying $i \neq j$.
The continuity of $f$ implies
$f$ ``pulls" $x_0$ into $S$, 
i.e., $f(x_0) \in S \setminus \{x_0\}$,
which, we saw above, is impossible.
\item $f(y)=H(y,1) \in Y \setminus \{y\}$ for some $y \in Y \setminus \{x_0\}$. 
Then $G_1$ is a homotopy from
$\id_Y$ to a nonidentity function
on $Y$; this is impossible, since
$Y$ is rigid.
\end{itemize}
The hypotheses of the cases 
listed above exhaust all
possibilities. Since each case
yields a contradiction, we
must have $f = \id_X$. Thus
$(X,\kappa)$ is rigid.
\end{proof}

\section{Cold and freezing sets}
Let $(X,\kappa)$ be a digital
image. Let $n \in \N^*$. 
We say
$f \in C(X,\kappa)$ is an
{\em $n$-map} \cite{BxLtd}
if $x \in X$ implies there
is a $\kappa$-path in X of
length at most $n$ from
$x$ to $f(x)$.

The following was observed in the
proof of Proposition~2.20
of~\cite{BxLtd}.

\begin{lem}
\label{1implies0}
    Let $(X,\kappa)$
    be a digital image. Let
    $f \in C(X.\kappa)$ be
    a 1-map. Then
    $f$ is $\kappa$-homotopic to
    $\id_X$.
\end{lem}

\begin{prop}
    \label{rigd1is0}
    {\rm \cite{BxLtd}}
    Let $(X,\kappa)$ be
    a connected rigid 
    digital image. Then the
    only 1-map in 
    $C(X,\kappa)$ is $\id_X$.
\end{prop}

\begin{thm}
{\rm \cite{BxFPSets2}}
\label{coldIffFreezing}
    Let $(X,\kappa)$ be
    a connected rigid 
    digital image.
    Then $A \subset X$ is a freezing
    set for $(X,\kappa)$
    if and only if
    $A$ is a cold
    set for $(X,\kappa)$.
\end{thm}

The converse of
Theorem~\ref{coldIffFreezing} is not generally valid,
as the following shows.

\begin{exl}
    Let $X = [0,2]_{\Z}$.
    Then $(X,c_1)$ is not
    rigid. However, each cold
    set for $(X,c_1)$ is
    freezing.
\end{exl}

\begin{proof}
    It is easily seen that
    $(X,c_1)$ is not rigid.
    It is easily seen that
    $A_1 = \{0,2\}$ and $X$
    are cold sets that are
    freezing. We show there
    are no other cold sets
    by showing $A_1$ is
    contained in any cold set
    $A$ for $(X,c_1)$.
    
    Suppose $0 \not \in A$.
    Then the function
    \[ f(x)= \left \{
    \begin{array}{ll}
       2 & \mbox{if } x=0; \\
       x  & \mbox{if } x \neq 0,
    \end{array} \right .
    \]
    satisfies $f \in C(X,c_1)$,
    $f|_A=\id_A$, and 
    $0 \not \adjeq_{c_1} f(0)$. Thus
    $A$ is not cold.

    Similarly, if
    $2 \not \in A$ then
    $A$ is not cold.
    Thus $A_1 \subset A$.
\end{proof}

\begin{thm}
    Let $(X,\kappa)$
    be a digital image. 
    Then $X$ is rigid if and
    only if the only 1-map in
    $C(X,\kappa)$ is $\id_X$.
\end{thm}

\begin{proof}
    If $X$ is rigid, it
    follows from
    Lemma~\ref{1implies0}
    that the only 1-map in
    $C(X,\kappa)$ is 
    $\id_X$.

Suppose the only 1-map in
    $C(X,\kappa)$ is 
    $\id_X$. Let 
$H: X \times [0,m]_{\Z} \to X$ be a
homotopy from $\id_X$
to $g \in C(X,\kappa)$.
We argue by induction
to show each induced
map $H_t(x)=H(x,t)$ 
is $\id_X$.

Clearly $H_0 = \id_X$.
Suppose $H_k=\id_X$ for
some $k$, $0 \le k < m$.
Then the continuity
properties of the homotopy 
$H$ imply $H_{k+1}$ is a 
1-map. By Proposition~\ref{rigd1is0}, 
$H_{k+1}=\id_X$. This
completes the induction.

Hence $g=H_m=\id_X$.
This shows $X$ is rigid.
\end{proof}

\section{Further remarks}
We have
studied implications of
the properties of
irreducibility and 
rigidity in digital
images for
Cartesian products,
wedges, and cold and
freezing sets.

\bibliographystyle{amsplain}

\end{document}